\documentclass[12pt]{amsart}
\usepackage{amssymb,amsmath,amsthm}
\newtheorem{thm}{Theorem}[section]
\newtheorem{lem}[thm]{Lemma}

\newtheorem{prop}[thm]{Proposition}

\newtheorem{defn}[thm]{Definition}

\newcommand{\kp}{k^\varphi}

\begin{document}
\title[Schur-Agler mappings of the unit ball]{An improved Julia-Caratheodory theorem for Schur-Agler mappings of the unit ball}
\author{Michael~T.~Jury}
\address{Department of Mathematics,
        University of Florida, 
        Gainesville, Florida 32603}
\email{mjury@math.ufl.edu}
\thanks{Partially supported by NSF grant DMS-0701268}
\date{\today}
\begin{abstract}We adapt Sarason's proof of the Julia-Caratheodory theorem to the class of Schur-Agler mappings of the unit ball, obtaining a strengthened form of this theorem.  In particular those quantities which appear in the classical theorem and depend only on the component of the mapping in the complex normal direction have $K$-limits (not just restricted $K$-limits) at the boundary.  
\end{abstract}
\maketitle

Let $\mathbb B^n$ denote the open unit ball in $n$-dimensional complex space.  In this note we show that holomorphic mappings $\varphi:\mathbb B^n\to \mathbb B^m$ belonging to the \emph{Schur-Agler class} (defined below) satisfy a strengthened form of the Julia-Caratheodory theorem (Theorem~\ref{T:stronger_jc}).  While the Schur-Agler class has received much attention in the past several years from operator theorists, relatively little seems to be known about the function-theoretic behavior of this class.    

For many operator theoretic applications, the Schur-Agler classes $\mathcal S(n,1)$ and $\mathcal S(n,n)$ are more appropriate analogues of the unit ball of $H^\infty (\mathbb{D})$ than are the larger classes $\text{Hol}(\mathbb B^n, \mathbb D)$ and $\text{Hol}(\mathbb B^n, \mathbb B^n)$.  For example, the Schur-Agler class is a natural setting for multivariable versions of von Neumann's inequality \cite{MR717226}, the Sz.-Nagy dilation theorem \cite{MR1668582}, commutant lifting theorems \cite{MR1846055} and the Nevanlinna-Pick interpolation theorem \cite{MR1882259}.  Additionally, every self-map of the ball belonging to the Schur-Agler class induces a bounded composition operator on the standard holomorphic function spaces \cite{jury-pams}, which is not true of general self-maps of the ball.  This last fact suggests that mappings in the Schur-Agler class should also enjoy function-theoretic privileges over generic maps of the ball, and is the motivation for this paper.  

Indeed there seems to be little known about the function theory of $\mathcal S(n,m)$ apart from what is true generically.  Recently Anderson, Dritschel and Rovnyak \cite{anderson-dritschel-rovnyak} have established a family of inequalities for derivatives of Schur-Agler functions, though it is not known if these inequalities hold generically.  In this paper we show that the Schur-Agler class satisfies a form of the Julia-Caratheodory theorem that is strictly stronger than what is true for general holomorphic functions on the unit ball.  The result is proved by adapting Sarason's Hilbert space proof of the classical Julia-Caratheodory theorem \cite[Chapter VI]{MR1289670} to the ball.  In fact Sarason's proof cannot prove the general Julia-Caratheodory theorem in higher dimensions, since it exploits the positivity of the de Branges-Rovnyak kernel.  The analogous kernel in several variables need no longer be positive, but since the Schur-Agler class is precisely the class for which this kernel is positive, the proof goes through but in fact proves a stronger result.

\begin{defn}
The \emph{Schur-Agler class} $\mathcal{S}(n,m)$ is the set of all holomorphic mappings $\varphi:\mathbb{B}^n \to \mathbb{B}^m$ such that the Hermitian kernel
\begin{equation*}
k^\varphi(z,w)=\frac{1-\langle \varphi(z), \varphi(w)\rangle}{1-\langle z,w \rangle}
\end{equation*}
is positive semidefinite.
\end{defn}
The kernel $k^\varphi$ is called the \emph{de~Branges-Rovnyak kernel} associated to $\varphi$.
We let $H^2_n$ denote the Hilbert space of holomorphic functions on $\mathbb{B}^n$ with reproducing kernel 
$$
k(z,w)=\frac{1}{1-\langle z,w \rangle}
$$
When $n>1$ the space $H^2_n$ is strictly smaller than the classical Hardy space $H^2$ (defined by spherical means); however in many ways it is the higher-variable analogue of $H^2(\mathbb D)$ appropriate for multivariable operator theory, see e.g. \cite{MR1882259, MR1668582, MR1846055}.  In this context, as mentioned above, the Schur-Agler classes play the role of the unit ball of the algebra of bounded analytic functions in $\mathbb{D}$, though we stress that when $n>1$ the inclusion $\mathcal{S}(n,m)\subset \text{Hol}(\mathbb B^n, \mathbb B^m)$ is always proper.  

Given a Schur-Agler mapping $\varphi \in \mathcal{S}(n,m)$, we can define another Hilbert function space $\mathcal{H}(\varphi)$ to be the space of holomorphic functions on $\mathbb{B}^n$ with reproducing kernel $k^\varphi$.  This space is always contractively contained in $H^2_n$:
\begin{lem}
If $\varphi \in S(n,m)$ and $f\in\mathcal{H}(\varphi)$ then $f\in H^2_n$ and 
$$
\|f\|_{H^2_n} \leq \|f\|_{\mathcal{H}(\varphi)}
$$
\end{lem}
\begin{proof}
The positivity of $k^\varphi$ implies that the operator 
$$
T:(f_1, \dots f_m) \to \sum_{k=1}^m \varphi_k f_k
$$
is contractive from the direct sum of $m$ copies of $H^2_n$ to $H^2_n$.  The de~Branges-Rovnyak kernel may then be written as 
$$
k^\varphi(z,w)=\big\langle(I-TT^*)^{1/2}k_w,(I-TT^*)^{1/2}k_z  \big\rangle_{H^2_n}  
$$
Now let $f\in \mathcal H(\varphi)$.  It follows from the standard de Branges-Rovnyak construction applied to $T$ \cite[Chapter 1]{MR1289670} that there exists $g\in H^2_n$ such that $f= (I-TT^*)^{1/2}g$ and $\|f\|_{\mathcal H(\varphi)}=\|g\|_{H^2_n}$.  Thus $f\in H^2_n$ and $\|f\|_{\mathcal H(\varphi)}\geq \|f\|_{H^2_n}$.
\end{proof}
We will be examining the boundary behavior of Schur-Agler mappings and to a lesser extent the behavior of functions in $\mathcal{H}(\varphi)$.  We recall here some basic notions in the study of boundary behavior of holomorphic functions on the unit ball, and refer to Rudin \cite[Chapter 8]{MR601594} (or Krantz \cite[Section 8.6]{MR1846625}) for details.

Given a point $\zeta\in\partial \mathbb{B}^n$ and a real number $\alpha >0$, the \emph{Koranyi region} $D_\alpha(\zeta)$ is the set
$$
D_\alpha(\zeta) =\{ z\in\mathbb{B}^n : |1-\langle z,\zeta\rangle |\leq \frac{\alpha}{2} (1-|z|^2)\}
$$
A function $f:\mathbb{B}^n\to \mathbb{C}$ has \emph{$K$-limit} $L$ at $\zeta$ if $\lim_{z\to \zeta}f(z)=L$ whenever $z$ tends to $\zeta$ within a Koranyi region.   Note that when $n=1$, a $K$-limit is just a nontangential limit; however for $n>1$ $K$-limits allow for parabolic approach in directions orthogonal to $\zeta$.  We shall also require the notion of a \emph{restricted $K$-limit}:  to define this, fix a point $\zeta\in\partial\mathbb{B}^n$ and consider a curve $\Gamma:[0, 1)\to \mathbb{B}^n$ such that $\Gamma(t)\to \zeta$ as $t\to 1$.  Let $\gamma(t)=\langle \Gamma(t), \zeta\rangle\zeta$ be the projection of $\Gamma$ onto the complex line through $\zeta$.  The curve $\Gamma$ is called \emph{special} if 
\begin{equation}
\lim_{t\to 1}\frac{|\Gamma -\gamma|^2}{1-|\gamma|^2}=0
\end{equation}
and \emph{restricted} if it is special and in addition
\begin{equation}
\frac{|\zeta-\gamma|}{1-|\gamma|^2}\leq A
\end{equation}
for some constant $A>0$.  We say that a function $f:\mathbb{B}^n\to \mathbb{C}$ has \emph{restricted $K$-limit} $L$ at $\zeta$ if $\lim_{z\to \zeta}f(z)=L$ along every restricted curve.  
\begin{lem}
If $f\in H^2_n$ then
$$
|f(z)|=o((1-|z|^2)^{-1/2})
$$
as $|z|\to 1$. 
\end{lem}
\begin{proof}
The Hilbert space norm of the reproducing kernel $k_z$ is 
$$
\sqrt{k(z,z)}=(1-|z|^2)^{-1/2}.
$$
The statement of the lemma is thus equivalent to the statement that the normalized kernel functions $\tilde {k_z}=k_z/\|k_z\|$ tend weakly to $0$ as $|z|\to 1$.  That this is the case follows readily from two observations:  1) if $f\in H^2_n$ is bounded, then $\langle f, \tilde{k_z}\rangle \to 0$ since $\|k_z\|\to \infty$, and 2) the bounded functions belonging to $H^2_n$ (e.g. the polynomials) are norm dense in $H^2_n$.  
\end{proof}
\begin{prop}
Suppose $\varphi\in\mathcal S(n,m)$ and $\zeta\in\partial\mathbb{B}^n$.   If 
$$
h(z)=\frac{1-\langle \varphi(z),\xi\rangle}{1-\langle z, \zeta\rangle}
$$
belongs to $\mathcal{H}(\varphi)$ for some $\xi\in \mathbb{C}^m$, then $|\xi|=1$ and $\varphi$ has $K$-limit $\xi$ at $\zeta$.
\end{prop}
\begin{proof}
If $h\in\mathcal{H}(\varphi)$ then by growth lemma $|h(z)|=o((1-|z|^2)^{-1/2}$.  So 
$$
|1-\langle \varphi(z),\xi\rangle| =o\left(\frac{|1-\langle z, \zeta\rangle|}{1-|z|^2}(1-|z|^2)^{1/2}\right)
$$
which goes to $0$ as $z\to \zeta$ within a Koranyi region; this establishes the claim.
\end{proof}
We are interested in Schur-Agler mappings satisfying the following condition, which we call condition (C) following Sarason:

\begin{equation*}
\tag{C}L=\liminf_{z\to \zeta}\frac{1-|\varphi(z)|^2}{1-|z|^2} < \infty 
\end{equation*}

The following is then the analogue, for the Schur-Agler class, of Sarason's Hilbert space formulation of the Julia-Caratheodory theorem \cite[Theorem VI-4]{MR1289670}:
\begin{thm}\label{T:sarason}
Let $\varphi\in\mathcal S(n,m)$ and $\zeta\in\partial\mathbb B^n$.  Then the following are equivalent:
\begin{enumerate}
\item Condition (C).

\item There exists $\xi\in\partial\mathbb{B}^m$ such that the function
$$
h(z)= \frac{1-\langle \varphi(z),\xi\rangle}{1-\langle z, \zeta\rangle}
$$
belongs to $\mathcal{H}(\varphi)$.
\item Every $f\in\mathcal{H}(\varphi)$ has a finite $K$-limit at $\zeta$.
\end{enumerate}
\end{thm}
\begin{proof}
First, suppose condition (C) holds.  Then there exists a sequence $z_n\to \zeta$ such that 
$$
L=\lim \|\kp_{z_n}\|^2_\varphi
$$
and by passing to a subsequence we may assume that $\varphi(z_n)\to \xi$ for some $\xi$ (necessarily $|\xi|=1$). By weak compactness of the closed unit ball in $\mathcal{H}(\varphi)$ (passing to a further subsequence if necessary) we have $\kp_{z_n}\to h$ weakly for some $h\in\mathcal{H}(\varphi)$.  Thus for all $z\in\mathbb{B}^n$, 
\begin{align*}
h(z)=\langle h, \kp_z\rangle_\varphi &= \lim_{n\to \infty} \langle \kp_{z_n},\kp_z\rangle_\varphi \\
&= \lim_{n\to \infty} \frac{1-\langle \varphi(z),\varphi(z_n)\rangle}{1-\langle z, z_n \rangle} \\
&=  \frac{1-\langle \varphi(z),\xi\rangle}{1-\langle z, \zeta\rangle}
\end{align*}
which proves (2).

Now assume (2).  By the lemma, $\varphi$ has $K$-limit $\xi$ at $\zeta$; we will write $\alpha=\varphi(\zeta)$ and $\kp_\zeta$ for the function $h$ in (2).  To prove (3) it suffices to prove that $\kp_z\to \kp_\zeta$ weakly as $z\to \zeta$ within a Koranyi region.  By taking inner products with the kernel functions $\kp_w$ it is clear that $\kp_z\to \kp_\zeta$ pointwise on $\mathbb{B}^n$ as $z\to \zeta$ in a Koranyi region.  Since the kernel functions $\kp_w$ span $\mathcal{H}(\varphi)$, it suffices to prove that the norms $\|\kp_z\|_\varphi$ remain bounded as $z\to \zeta$ in a Koranyi region.  For each $z\in\mathbb{B}^n$ we have
$$
\langle \kp_\zeta ,\kp_z\rangle =\frac{1-\langle \varphi(z),\varphi(\zeta)\rangle}{1-\langle z, \zeta\rangle}
$$
so by the Cauchy-Schwarz inequality
$$
\left| \frac{1-\langle \varphi(z),\varphi(\zeta)\rangle}{1-\langle z, \zeta\rangle}\right|^2 \leq \|\kp_\zeta\|_\varphi^2 \|\kp_z\|^2_\varphi = \|\kp_\zeta\|_\varphi^2 \left(\frac{1-|\varphi(z)|^2}{1-|z|^2} \right)
$$
The numerator on the left hand side dominates $(1-|\varphi(z)|)^2$, so 
$$
\frac{(1-|\varphi(z)|)^2}{|1-\langle z, \zeta\rangle|^2} \leq \|\kp_\zeta\|_\varphi^2 \left(\frac{1-|\varphi(z)|^2}{1-|z|^2} \right)
$$
which implies
$$
\|\kp_z\|_\varphi^2 = \frac{1-|\varphi(z)|^2}{1-|z|^2} \leq \|\kp_\zeta\|_\varphi^2 \left( \frac{1+|\varphi(z)|^2}{1+|z|^2}\right)\left( \frac{|1-\langle z, \zeta\rangle|}{1-|z|}\right)^2
$$
The right hand side remains bounded as $z\to \zeta$ in a Koranyi region, which proves (3).

The proof that (3) implies (1) is immediate, since by the principle of uniform boundedness the norms $\|\kp_z\|_\varphi$ stay bounded as $z\to \zeta$ in a Koranyi region, which implies condition (C).
\end{proof}

\begin{thm}\label{T:sarason2}
Suppose $\varphi\in\mathcal S(n,m)$ and satisfies condition (C).  Then the function
$$
\frac{1-|\langle \varphi(z),\xi\rangle|^2}{1-|\langle z,\zeta\rangle|^2}
$$
has $K$-limit $L$ at $\zeta$.  
\end{thm}
\begin{proof}
By pre- and post-composing with unitary rotations, we may assume without loss of generality that (in the nomenclature of previous theorem) $\xi=e_1$ and $\zeta=e_1$.  (We are using $e_1$ to refer to vectors in two different spaces, but this should cause no confusion.)  

Starting with the identity
$$
1-\varphi_1(z) =(1-z_1)\langle \kp_{e_1}, \kp_z \rangle
$$
we find
$$
|\varphi_1(z)|^2 = 1-2\text{Re}[(1-z_1) \langle \kp_{e_1}, \kp_z \rangle] +|1-z_1|^2 |\langle \kp_{e_1}, \kp_z \rangle |^2
$$
From what has already been proved, the last term is $o(1-|z_1|^2)$ as $z\to e_1$ within a Koranyi region.  Thus
$$
K\text{-}\lim_{z\to e_1}\frac{1-|\varphi_1(z)|^2}{1-|z_1|^2}=K\text{-}\lim_{z\to e_1} \frac{2\text{Re}[(1-z_1) \langle \kp_{e_1}, \kp_z \rangle]}{1-|z_1|^2}
$$
As $z\to e_1$ in a Koranyi region, the real part of 
$$
\frac{1-z_1}{1-|z_1|^2}
$$
tends to $1/2$ and its imaginary part remains bounded.  The real part of $\langle \kp_{e_1}, \kp_z \rangle $ tends to $\|\kp_{e_1}\|^2$ and its imaginary part tends to $0$.  Thus 
$$
K\text{-}\lim_{z\to e_1} \frac{2\text{Re}[(1-z_1) \langle \kp_{e_1}, \kp_z \rangle]}{1-|z_1|^2} = \|\kp_{e_1}\|^2=L
$$
which completes the proof.
\end{proof}

Combining statements (2) and (3) of Theorem~\ref{T:sarason} we obtain our first strengthened conclusion, namely that the function $h$ has finite $K$-limit at $\zeta$.  For general $\varphi$ this will exist only as a restricted $K$-limit.  The same is true for the expression in Theorem~\ref{T:sarason2}.  These facts will allow us to strengthen the convergence results for directional derivatives of the component of $\varphi$ in the $\zeta$ direction.

In the disk, Theorem~\ref{T:sarason2} says that $\|\kp_z\|_\varphi \to \|\kp_\zeta\|_\varphi$ as $z\to \zeta$ nontangentially; together with the weak convergence of $\kp_z$ to $\kp_\zeta$ this shows that in fact $\kp_z \to \kp_\zeta$ in norm.  In the ball we would like to establish $\|\kp_z\|_\varphi \to \|\kp_\zeta\|_\varphi$ or equivalently
$$
\frac{1-|\varphi(z)|^2}{1-|z|^2}\to L
$$
in as general a sense as possible.  For generic self-maps $\varphi$, this limit exists restrictedly but not as a $K$-limit in general.  Unlike the previous results, however, this cannot be improved for Schur-Agler mappings; in fact for the Schur-Agler mapping $\varphi(z)=z_1$ the above expression does not have a $K$-limit at $e_1$.  Thus in the ball we only have $\|\kp_z\|_\varphi \to \|\kp_\zeta\|_\varphi$ (and hence $\kp_z \to \kp_\zeta$ in norm) when $z\to \zeta$ restrictedly.

The following is Rudin's version of the Caratheodory theorem on the ball:

\begin{thm}\label{T:rudin}
Suppose $\varphi=(\varphi_1, \dots \varphi_m)$ is a holomorphic mapping from $\mathbb{B}^n$ to $\mathbb{B}^m$ satisfying condition (C)at $e_1$.  Suppose $2\leq j\leq m$ and $2\leq k\leq n$.  The following functions are then bounded in every Koranyi region $D_\alpha(e_1)$:
\begin{itemize}
\item[(i)]  $(1-\varphi_1(z))/(1-z_1)$
\item[(ii)]  $(D_1 \varphi_1)(z)$
\item[(iii)]  $\varphi_j(z)/(1-z_1)^{1/2}$
\item[(iv)]  $(1-z_1)^{1/2}(D_1 \varphi_j)(z)$
\item[(v)]  $(D_k\varphi_1)(z)/(1-z_1)^{1/2}$
\item[(vi)]  $(D_k \varphi_j)(z)$
\end{itemize}
Moreover, the functions (i), (ii) have restricted K-limit $L$ at $e_1$, and the functions (iii), (iv), (v) have restricted K-limit $0$ at $e_1$.
\end{thm}

We next show that for $\varphi\in\mathcal S(n,m)$, the restricted $K$-limits in (i), (ii) and (v) can be improved to $K$-limits.  Note that these are precisely the expressions that involve only the $e_1$ component of $\varphi$.  This is to be expected, since the improvement derives from the fact that the kernel $k^\varphi_\zeta$ has a $K$-limit at $\zeta$, and this kernel depends only on the component of $\varphi$ in the $\zeta$ (that is, the complex normal) direction.  Indeed, the limits of (iii) and (iv) cannot be improved to K-limits, since the counterexamples given in \cite{MR601594} are in fact Schur-Agler mappings; this will be shown after proving the next theorem.  Before beginning we recall Lemma 8.5.5 of \cite{MR601594} which will be used in the proof.

\begin{lem}
Suppose $1<\alpha <\beta$, $\delta=\frac13 (1/\alpha -1/\beta)$, and $z=(z_1, z^\prime)\in D_\alpha$.  
\begin{itemize}
\item[(i)]  If $|\lambda|\leq \delta|1-z_1|$ then $(z_1 +\lambda, z^\prime)\in D_\beta$.
\item[(ii)]  If $|w|\leq \delta|1-z_1|^{1/2}$ then $(z_1, z^\prime +w^\prime)\in D_\beta$.
\end{itemize}
\end{lem}

\begin{thm}\label{T:stronger_jc}
Suppose that $\varphi\in\mathcal S(n,m)$ and satisfies condition (C).  Then in (i), (ii) and (v) of Theorem~\ref{T:rudin}, restricted K-limit can be improved to K-limit.
\end{thm}

\begin{proof}
Since we are assuming condition (C), we know from statement (2) of Theorem~\ref{T:sarason} that the function
$$
k_{e_1}^\varphi (z)=\frac{1-\varphi_1(z)}{1-z_1}
$$
belongs to $\mathcal{H}(\varphi)$ and hence by statement (3) has a $K$-limit at $e_1$; this limit must of course equal $L$.  

For (ii), suppose $1<\alpha<\beta$, choose $\delta$ as in the lemma, let $z\in D_\alpha$ and put
$$
r=r(z)=\delta |1-z_1|
$$
As in \cite{MR601594}, express $D_1\varphi_1$ using the Cauchy formula; after some manipulation we obtain
\begin{equation}\label{E:D1}
(D_1 \varphi_1)(z)=\frac{1}{2\pi}\int_{-\pi}^\pi \frac{1-\varphi_1(z_1 +re^{i\theta}, z^\prime)}{1-(z_1+re^{i\theta})}\cdot \left\{1-\frac{1-z_1}{re^{i\theta}}  \right\}\, d\theta
\end{equation}
We must show that the above expression tends to $L$ along any sequence converging to $e_1$ within $D_\alpha$; in fact it suffices to show that given any such sequence, $D_1\varphi_1$ converges to $L$ along some subsequence.  In particular we may assume that we have chosen a sequence $(z_n)$ such that
$$
\lim_{n\to \infty} \frac{1-z_{n,1}}{r(z_n)e^{i\theta}}=\frac{1}{\delta e^{i\theta}}\lim_{n\to \infty}\frac{1-z_{n,1}}{|1-z_{n,1}|} 
$$
exists, and is equal to some complex number $\lambda$.  Then as $z_n\to e_1$ in $D_\alpha$, we have 
$$
z_n +r(z_n)e^{i\theta} e_1 \to e_1
$$
in $D_\beta$, so the integrand in (\ref{E:D1}) converges to 
$$
L\cdot (1-\frac{\lambda}{\delta e^{i\theta}})
$$
for every $\theta\in [0, 2\pi]$.  Since this integrates to $L$, and the integrands are uniformly bounded, we conclude $D_1\varphi_1(z_n)\to L$ by the dominated convergence theorem.  

The $K$-limit of (v) is established similarly:  we let $\alpha, \beta, \delta$ be as before, and for $z\in D_\alpha(e_1)$ we define
$$
\rho=\rho(z)=\delta |1-z_1|^{1/2}
$$
Then by the lemma, $(z_1 ,z^\prime +w^\prime)\in D_\beta(e_1)$ for all $w^\prime$ with $|w^\prime|\leq \rho$.  Assuming $k=2$ (without loss of generality), we apply the Cauchy formula to obtain for every $z\in D_\alpha$
\begin{equation*}
\frac{(D_2 \varphi_1)(z)}{(1-z_1)^{1/2}} =-\frac{(1-z_1)^{1/2}}{\rho(z)}\cdot \frac{1}{2\pi} \int_{-\pi}^\pi \frac{1-\varphi_1(z_1,z_2+\rho e^{i\theta},\dots)}{1-z_1}e^{i\theta}\, d\theta
\end{equation*} 
The factor outside the integral is bounded.  As $z\to e_1$ within $D_\alpha$, $z+\rho(z)e^{i\theta}e_2 \to e_1$ within $D_\beta$, so the integrand tends to $Le^{i\theta}$ for every $\theta$.  Thus  
$$
\frac{(D_2 \varphi_1)(z)}{(1-z_1)^{1/2}} \to 0
$$
by the dominated convergence theorem.  
\end{proof}

Rudin \cite{MR601594} gives counterexamples to show that ``restricted $K$-limit'' cannot be improved to ``$K$-limit'' in Theorem~\ref{T:rudin}; in the case of (iii) and (iv), the example is a map $\varphi:\mathbb{B}^2\to \mathbb{B}^2$ of the form
\begin{equation}\label{E:example}
\varphi (z_1, z_2) =(z_1, z_2 g(z_1))
\end{equation}
for a suitably chosen holomorphic function $g:\mathbb{D}\to \mathbb{D}$.  It is not hard to show that any map of the form (\ref{E:example}) belongs to $\mathcal{S}(2,2)$.  To see this, first observe that because $g:\mathbb{D}\to \mathbb{D}$, the kernel
$$
\frac{1-g(z)\overline{g(w)}}{1-z\overline{w}}
$$
is positive.  We may then write
\begin{align}
\frac{1-\langle \varphi(z), \varphi(w)\rangle }{1-\langle z,w\rangle} &= \frac{1-\langle z,w\rangle +z_2 \overline{w_2}-z_2 \overline{w_2}g(z_1)\overline{g(w_1)}}{1-\langle z,w\rangle} \\
&= 1+ z_2\overline{w_2} \frac{1-g(z_1)\overline{g(w_1)}}{1-z_1\overline{w_1}}\cdot \frac{1-z_1\overline{w_1}}{1-\langle z,w\rangle}
\end{align}
which is positive.

\bibliographystyle{plain} 
\bibliography{jc} 
\end{document}